\documentclass[twocolumn,10pt]{asme2ej}
\usepackage{ifpdf}
\usepackage{multirow}
\ifpdf
    \usepackage{graphicx}
    \usepackage{epstopdf}
\else
    \usepackage{graphicx}

\fi

\usepackage{mathptmx} % assumes new font selection scheme installed
\usepackage{times} % assumes new font selection scheme installed
\usepackage{amsmath} % assumes amsmath package installed
\usepackage{amssymb}  % assumes amsmath package installed
\usepackage{tabularx}
\usepackage{algorithm}
\usepackage{algorithmic}
\usepackage{ulem}
\usepackage{dsfont}

\usepackage[noadjust]{cite}

\usepackage[dvipsnames,usenames]{color}
\usepackage[colorlinks,%
linkcolor=BrickRed,%
filecolor=BrickRed,%
citecolor=RoyalPurple,%
]{hyperref}
\usepackage{caption}
\usepackage{subfigure}
\newcounter{rmnum}

\usepackage{color}

\usepackage{tabularx}
\def\notes#1{\marginpar{\tiny #1}\typeout{Notes!
Notes!
Notes!
}}
\renewcommand{\notes}[1]{\typeout{notes!}}

\def\beq{\begin{eqnarray}} 
\def\bc{\begin{center}} 
\def\be{\begin{enumerate}}
\def\bi{\begin{itemize}} 
\def\bs{\begin{small}}
\def\bS{\begin{slide}}
\def\ec{\end{center}} 
\def\ee{\end{enumerate}}
\def\ei{\end{itemize}}
\def\es{\end{small}}
\def\eS{\end{slide}}
\def\eeq{\end{eqnarray}}

\title{Interacting particle filters for simultaneous state and parameter estimation}

\author{Angwenyi David, 
    \affiliation{
Institut f\"ur Mathematik \\ 
    Universit\"at Potsdam\\
Email: kipkoej@gmail.com
    }	
}

\author{Jana de Wiljes,
    \affiliation{
    Institut f\"ur Mathematik \\ 
    Universit\"at Potsdam\\
    Email: wiljes@uni-potsdam.de
    }	
}

\author{Sebastian Reich, 
    \affiliation{
	Institut f\"ur Mathematik\\ 
	Universit\"at Potsdam and\\
	Department of Mathematics and Statistics \\
	University of Reading\\
	Email: sereich@uni-potsdam.de
    }	
}
\begin{document}
\maketitle

\vspace{-0.5in}
\begin{abstract} Simultaneous state and parameter estimation arises from various 
applicational areas but presents a major computational challenge. Most available 
Markov chain or sequential Monte Carlo techniques are applicable to relatively 
low dimensional problems only. Alternative methods, such as the ensemble Kalman filter or other
ensemble transform filters have, on the other hand, been successfully applied to high dimensional
state estimation problems. In this paper, we propose an extension of these techniques to 
high dimensional state space models which depend on a few unknown parameters. More
specifically, we combine the ensemble Kalman-Bucy filter for the continuous-time filtering problem
with a generalized ensemble transform particle filter for intermittent parameter updates.
We demonstrate the performance of this two stage update filter for a wave equation with unknown wave velocity parameter.\end{abstract}

\section{Introduction}
%
%%%%%%%%%%%%%%%%%%%%%%%%%%%%%%%%%%%%%%%%%%%%%%%%%%%%
There is a high demand across different disciplines for methods that allow for efficient and reliable state-parameter estimation for high-dimensional and nonlinear evolution equations.
While the theoretical foundation of state and parameter estimation for stochastic differential equations (SDEs) is well established (see, for example, \cite{sr:jazwinski,sr:crisan} and efficient computational
methods for low dimensional problems are available (see, for example, \cite{sr:Doucet, KDSMC15}), joint state-parameter estimation for high dimensional problems remains an area of active
research. A major breakthrough in that direction has been achieved through the development of the ensemble Kalman filter (EnKF) for state estimation of discretized partial differential equation
models arising, for example, from meteorology and oil reservoir exploration \cite{sr:evensen,sr:stuart15,Reich2015}. The success of the EnKF has triggered the development of 
a large variety of related ensemble transform filters with the aim of removing the underlying Gaussian distributional assumptions of the EnKF. Here we wish to mention in particular the work
of \cite{sr:crisan10,sr:crisan11,sr:br11,sr:meyn13,jdw:TaghvaeideWiljesMehtaReich} on the continuous-time filtering problem and \cite{Xiong2006,sr:reich13,Toedter2015,jdw:DeWiljesAcevedoReich2017}
on the intermittent filtering problem. 

In this paper, we propose an extension of the ensemble transform filtering approach to the continuous-time combined state and parameter estimation problem. Instead of applying an ensemble transform
filtering approach directly to the extended state-parameter phase space, we propose to exploit the particular structure of the joint conditional distribution and approximate it via a hybrid {\it ansatz} 
combining two different interacting particle filters; namely the ensemble Kalman-Bucy filter (EnKBF) \cite{sr:br11} for state estimation and the ensemble transform particle filter (ETPF) \cite{sr:reich13,Reich2015} for the parameter updates. 
Such an approach is advantageous provided the distribution in the states given model parameters is nearly Gaussian while the distribution in the parameters may be non-Gaussian. Furthermore, the 
main additional computational complexity arises from the update of the parameters through an appropriate extension of the ETPF. Here we assume that the number of unknown parameters is much smaller
than the dimension of state space of the underlying SDE model.

The remainder of the paper is structured as follows. In Section \ref{sec:PROBLEM}, we will discuss the theoretical foundation of the considered Bayesian inference problem 
and formulate the basic algorithmic approach. Our proposed approach for a sequential update of the model parameters and the required extension of the ETPF is provided in 
Section \ref{Sec:ETPF}. A summary of the overall algorithm is provided in Section \ref{sec:AS} and numerical results for a stochastic wave equation in Section \ref{sec:NE}.
Our conclusions can be found in Section \ref{Sec:Conclusions}.

%%%%%%%%%%%%%%%%%%%%%%%%%%%%%%%%%%%%%%%%%%%%%%
%
\section{Problem formulation and proposed ansatz} \label{sec:PROBLEM}
%
%%%%%%%%%%%%%%%%%%%%%%%%%%%%%%%%%%%%%%%%%%%%%%

We consider the following time-continuous filtering problem: estimate a reference trajectory $x_{\rm ref}(t) \in \mathbb{R}^{N_x}$ and  $t\in [0,T]$ and a vector of unknown reference parameters $\lambda_{\rm ref}\in\mathbb{R}^{N_{\lambda}}$ of a SDE
\begin{equation} \label{sde}
{\rm d} x_t = f(x_t,\lambda) {\rm d} t + Q^{1/2} {\rm d} W_t
\end{equation}
from continuous-time observations
\begin{equation} \label{obs}
{\rm d}y_t = h(x_t) {\rm d} t + R^{1/2} {\rm d} V_t.
\end{equation}
Here both $W_t \in \mathbb{R}^{N_x}$ and $V_t\in \mathbb{R}^{N_y}$ denote standard multi-dimensional Brownian motions.
A common approach to joined parameter and state estimation is to augment the SDE (\ref{sde}) by the trivial dynamics
\begin{align} \label{sde2}
{\rm d} \lambda_t=  0
\end{align}
in the parameters $\lambda$. Then the distribution of interest is the conditional density $\tilde\pi_t(z):=\pi_t(z|y_{[0,t]})$ in the augmented state variable 
$z=(x^{\top},\lambda^{\top})^{\top} \in\mathbb{R}^{N_z}$. The time evolution of the conditional density $\tilde \pi_t$ is described by the 
Kushner-Stratonovich equation \cite{sr:jazwinski,sr:crisan}, which we state in the form
\begin{align} \nonumber
\tilde\pi_t(g) = &\tilde\pi_0(g) + \int_0^t\tilde\pi_s(\mathcal{L}g){\rm d} s \\& \qquad +\int_0^t \left(\tilde\pi_s(gh)-\tilde\pi_s(g) \bar h_s \right)^{\rm T} R^{-1} 
\left({\rm d} y_s -\bar h_s {\rm d} s\right)\,,
\label{KZ}
\end{align} 
where $\bar h_s := \tilde\pi_s(h)$ and
\[
\mathcal{L}g := f \cdot \nabla_x g  + \frac{1}{2}\sum_{k,l=1}^{N_x} Q_{kl}\frac{\partial^2g}{\partial x_k \partial x_l}.
\]
Analytical solutions of (\ref{KZ}) are generally not available and sequential Monte Carlo (SMC) 
methods are often employed in order to approximate the marginal density 
by empirical measures. Here we follow recently developed SMC methods which rely on appropriately defined modified evolution equations for 
particles $z_t^l$, $l = 1,\ldots,L$, such that
\begin{align} \label{emp_marginal}
\tilde \pi_{t}(z) &\approx \frac{1}{N} \sum_{l=1}^N \delta (z-\tilde z_t^{l})\,.
\end{align}
The feedback particle filter (FPF) \cite{sr:meyn13,Chi_ACC16} 
is one of these, so called, particle flow filters, which 
is characterized by the modified SDE
 \begin{align} \label{FPF}
{\rm d} z^l_t &= \begin{bmatrix}f(z^l_t)dt + Q^{1/2} {\rm d}W^l_t\\ 0\end{bmatrix} +K^l_t\circ {\rm d}I_t^l
\end{align}
with gain factors $K_t^l = K_t(z_t^l)$ and innovation
\begin{equation} \label{innovation}
{\rm d}I_t^l = {\rm d}y_t-\frac{h(x^l_t)+\bar h_t}{2} {\rm d}t \,.
\end{equation}
Here the Stratonovitch interpretation of the SDE (\ref{FPF}) should be used \cite{pavliotis14}.
The gain function $K_t$ is determined by the elliptic partial differential equation
\begin{equation}\label{pdegainfactor}
-\nabla_z \cdot(\tilde\pi_{t}K_t)=  \tilde\pi_{t}R^{-1}(h-\bar h_t)^{\rm T}\,.
\end{equation}
Note that $\tilde\pi_{t}$ is unknown and needs to be approximated by (\ref{emp_marginal}). In other words, $K_t$ is typically found as a weak approximation
to (\ref{pdegainfactor}). Different numerical approaches for solving (\ref{pdegainfactor}) can be found in \cite{jdw:TaghvaeideWiljesMehtaReich}. We mention that the innovation (\ref{innovation}) can be
replaced by the alternative form
\begin{equation} \label{innovation2}
{\rm d}I_t^l = {\rm d}y_t- h(x^l_t) {\rm d}t + R^{1/2} {\rm d}U_t^l\,,
\end{equation}
where $U_t^l$ denote standard $N_y$-dimensional Brownian motion independent of 
$W_t$ and $V_t$. The statistical equivalence can be shown following the arguments of
Appendix A in \cite{jdw:TaghvaeideWiljesMehtaReich}.

While (\ref{FPF}) is very appealing, its numerical implementation can be demanding for high-dimensional systems which require a large number, $L$, of particles
$z_t^l$. In order to address this issue we propose to rewrite the joint distribution $\tilde \pi_t$ in its desintegrated form, {\rm i.e.}
\begin{equation} 
\tilde\pi_t(z)=\hat \pi_t(x|\lambda)\hat\pi_t(\lambda).
\end{equation}
A corresponding particle approximation can be defined as follows:
\begin{align} \label{mixture}
\tilde \pi_{t}(z)\approx \sum^L_{i=1}  w_t^i \delta (\lambda-\hat \lambda_0^{i})\frac{1}{M}\sum_{j=1}^M \delta (x-\hat x_t^{i,j})\,,
\end{align}
where $\hat \lambda_0^{i}\sim\hat\pi_0(\lambda)$, $i=1,\ldots,L$, are constant parameter values drawn from the prior parameter distribution with
time-dependent weights $w_t^i$, $i=1,\ldots,L$. There is also a set of $M$ time-dependent states $\{\hat x_t^{i,j}\}_{j=1}^M$ for each parameter vector
$\hat \lambda_0^i$, $i\in \{1,\ldots,L\}$. The evolution equations for these states are given by the 
FPF with the parameters $\hat \lambda_0^i$ held fixed, {\rm i.e.},
\begin{align}
{\rm d} \hat x^{i,j}_t &= f( \hat x^{i,j}_t,\hat \lambda^i_0) {\rm d} 
t + Q^{1/2} {\rm d}W_t^{i,j}  +\hat K^{i,j}_t\circ {\rm d}I_t^{i,j}\,,
\end{align}
where $\hat K^{i,j} = K_t(\hat x_t^{i,j}, \hat \lambda_0^i)$ is determined by an appropriate numerical approximation to
\begin{equation}
-\nabla_x \cdot(\hat\pi_{t}\hat K_t)= \hat\pi_{t}R^{-1}(h-\bar h^i_t)^{\rm T}
\end{equation}
and ${\rm d}I_t^{i,j}$ denotes the innovation, {\rm i.e.},
\begin{equation} \label{innovation1b}
{\rm d}I_t^{i,j} = {\rm d}y_t-\frac{h(\hat x^{i,j}_t)+\bar h^{i}_t}{2} {\rm d} t 
\end{equation}
or
\begin{equation} \label{innovation2b}
{\rm d}I_t^{i,j} = {\rm d}y_t- h(\hat x^{i,j}_t){\rm d} t + R^{1/2} {\rm d}U_t^{i,j} \,,
\end{equation}
respectively.
The time evolution of the normalized importance weights $w^i_t$ are calculated according to
\begin{equation}\label{eq:weights}
{\rm d} w_t^i =w_t^i(\bar h^i_{t}-\bar h_t)^{\rm T} R^{-1}({\rm d} y_t-\bar h_t {\rm d} t)
\end{equation}
with $w_0^i:=1/L$ initially and
\begin{align}
\bar h^{i}_t=\frac{1}{M}\sum_{j=1}^M \bar h(\hat x_t^{i,j})\,,\quad \quad \bar h_t=\frac{1}{L}\sum_{i=1}^L \bar h^{i}_t\, .
\end{align}
Note again that the parameter values are kept constant in (\ref{mixture}), {\rm i.e.}, $\hat 
\lambda_n^i = \hat \lambda_0^i$.

A special case of the FPF scheme arises when the gain factor is assumed to be constant, which results in the popular EnKBF:
\begin{equation}\label{ENKBF}
{\rm d} \hat x_t^{i,j}=f(\hat x_{t}^{i,j},\hat \lambda_0^i){\rm d} t + Q^{1/2} {\rm d} W_t^{i,j} + C^i_t R^{-1} {\rm d}I_t^{i,j}\,,
\end{equation}
where the covariance matrix $C_t^i$ is determined empirically, {\rm i.e.},
\[
C^i_{t}=\frac{1}{M-1}\sum_{j=1}^M(\tilde x^{i,j}_{t}-\bar x^{i}_{t})(h(\hat x^{i,j}_{t})-\bar h^{i}_{t})^{\top}\,,\,\,
\bar x^i_t=\frac{1}{M}\sum_{j=1}^M \hat x_t^{i,j}\,,
\]
and the innovation ${\rm d}I_t^{i,j}$ is either given by (\ref{innovation1b}) or (\ref{innovation2b}),
respectively.
The EnKBF produces asymptotically correct results in a linear model setting when the posterior is a Gaussian distribution but is also successfully employed for state estimation in the context of 
strongly nonlinear model scenarios \cite{sr:br11}. In this paper, we employ the EnKBF to forward state samples, $\hat x_t^{i,j}$, in time. In other words, we interpret (\ref{mixture}) as a 
weighted Gaussian mixture approximation to the conditional filtering distribution, $\tilde \pi_t$. 

The effective mixture size, defined by
\begin{equation} \label{effective_SS}
L_t^{\rm eff} = \frac{1}{\sum_{i=1}^L (w_t^i)^2}\,,
\end{equation}
will deteriorate as time progresses, in general. A classic approach would be to resample the parameter values $\hat \lambda_0^i$ jointly with their state samples $\{\hat x_t^{i,j}\}_{j=1}^M$, 
$i\in \{1,\ldots,L\}$ at an appropriate instance   $t=t^\ast$ of time according to their weights 
$w_{t^\ast}^i$ in order to produce  an equally weighted mixture (\ref{mixture}). However, 
resampling with replacement would produce identical sets of 
parameters and associated state samples. Hence, we propose an extension of the ETPF \cite{sr:reich13,Reich2015} to (\ref{mixture}). Contrary to the EnKBF, 
the ETPF produces a consistent approximation of the gain factor of the FPF on the basis of an 
optimal transport problem \cite{jdw:TaghvaeideWiljesMehtaReich}. 
The ETPF has also been shown to work well  under relatively small number of particles and high dimensional systems  when combined with localization 
\cite{jdw:DeWiljesAcevedoReich2017,Chustagulprom2015}.

%In other words, we will alternate between two numerical approximations of the FPF to produce simultaneous updates of the state and the parameter vector. 
%The next sections is devoted to introduce the key ideas of the ETPF and to briefly address its numerical implementation in the context of (\ref{mixture}). 

%%%%%%%%%%%%%%%%%%%%%%%%%%%%%%%%%%%%%
%
\section{Ensemble transform particle filter}\label{Sec:ETPF}
%
%%%%%%%%%%%%%%%%%%%%%%%%%%%%%%%%%%%%

As mentioned above, the ETPF is an numerical approximation of the feedback control law of the FPF induced by a linear transport problem \cite{jdw:TaghvaeideWiljesMehtaReich}. A different interpretation is that the ETPF  replaces the resampling step of the classical particle filter with a linear transformation \cite{Reich2015,sr:reich13}. The key idea is to choose a linear transformation that connects the 
empirical measure of the weighted prior ensemble with an equally weighted posterior ensemble in the sense of optimal transportation. Intuitively, one would like to achieve a high correlation between the prior and posterior samples. More generally, the optimal transport problem between two 
weighted empirical measures $\nu_1$ and $\nu_2$, given by
\[
\nu_1(y) = \sum_{i=1}^L w^i_1 \delta (y-y^i_1)\,,\quad \nu_2(y) = \sum_{i=1}^L w_2^i \delta (y-y^i_2),
\]
can be formulated as follows \cite{Reich2015}. Introduce the set
\[
U(W_1,W_2) = \{ T\in \mathbb{R}^{L\times L}: t_{ij}\ge 0,\,\sum_{j=1}^L t_{ij} w_1^i,\,\sum_{i=1}^L t_{ij} = w_2^j\}
\]
of admissible bi-stochastic matrices $T$ 
and the $L\times L$ matrix of mutually distances $M_{Y_1,Y_2}$ with entries
\[
(M_{Y_1,Y_2})_{ij} = \|y_1^i - y_2^j\|^2\,.
\]
Then the Wasserstein distance between $\nu_1$ and $\nu_2$ is defined by
\begin{equation} \label{wasserstein}
W_2^2(\nu_1,\nu_2) = \min_{T\in U(W_1,W_2)} {\rm tr}\,\left(T^{\rm T} M_{Y_1,Y_2} \right).
\end{equation}
The matrix $T^\ast \in U(W_1,W_2)$, which achieves the minimum in (\ref{wasserstein}), is called the optimal coupling between
$\nu_1$ and $\nu_2$. 

The ETPF relies on the special situation that the vector $W_1 = (w_1^1,\ldots,w_1^L)^{\rm T}$ represents the importance weights of the prior samples $y^i$, $i\in \{1,\ldots,L\}$, 
and $W_2 = (1/L,\ldots,1/L)^{\rm T}$. Furthermore, the ETPF also uses $Y_1 = Y_2 = Y := (y^1,\ldots,y^L)^{\rm T}$ and the desired equally weighted posterior samples are defined by
\[
\tilde y^j = M \sum_{i=1}^L  y^i t_{ij}^\ast \,.
\]
Solving an optimal transport problem is computationally demanding for large sample sizes $L$. 
This issue has been addressed in \cite{Cuturi2013,CuturiDoucet2013} 
via a Sinkhorn approximation which reduces the complexity of the optimal transport problem from $\mathcal{O}(L^3log(L))$ to $\mathcal{O}(L^2)$. It is shown in \cite{jdw:DeWiljesAcevedoReich2017}
how this approximation can be employed successfully in the context of sequential filtering.

In case of the weighted mixture approximation (\ref{mixture}), the ETPF is implemented at an appropriate instance $t=t^\ast$ of time as follow. 
First, we define the distance matrix $M_{Y_1,Y_2}$. There are two choices. Either
one sets $Y_1 = Y_2  = (\hat \lambda_0^1,\ldots, \hat \lambda_0^L)^{\rm T}$ or one uses the extended vectors $\hat z_t^i = ((\hat \lambda_0^i)^{\rm T},(\bar x_t^i)^{\rm T})^{\rm T} \in \mathbb{R}^{N_z}$
instead of $\hat \lambda_0^i$ in both $Y_1$ and $Y_2$. Second, the weight vector $W_1$ is defined by $W_1 = (w_t^1,\ldots,w_t^L)^{\rm T}$. 
Denoting the solution of the optimal transport problem again by $T^\ast$, equally weighted parameter values are finally provided by
\begin{equation} \label{transform_lambda}
\hat \lambda_t^j = L \sum_{i=1}^L \hat \lambda_0^i t^\ast_{ij}\,.
\end{equation}
We also need to transform the associated state samples $\hat x_t^{i,j}$. The obvious choice is
\begin{equation} \label{transform_x}
\hat x_{t^+}^{i,j} = M \sum_{l=1}^L \hat x_t^{l,j} t^\ast_{il} \qquad \forall j \in \{1,\ldots,M\}\,.
\end{equation}
This requires, however, that the state samples $\hat x_t^{i,j}$ are optimally correlated for each fixed index $i \in \{1,\ldots,L\}$. This can be achieved either through an appropriate initialization of
the state samples or through finding an appropriate permutation matrix 
$P^i \in \mathbb{R}^{M\times M}$ for each set of state samples $\{\hat x_t^{i,j}\}_{j=1}^M$ 
via an associated Wasserstein barycenter problem \cite{CuturiDoucet2013}. More specifically, introduce $L$ equally weighted empirical measures
\[
\nu_i(x) = \frac{1}{M}\sum_{j=1}^M \delta (x-\hat x_t^{i,j})
\]
and the empirical measure
\[
\nu (x) = \frac{1}{M} \sum_{i=1}^M \delta (x-x^j)
\]
with its locations $x^j$, $j=1,\ldots,M$, determined as the minimizer of the functional
\begin{equation}\label{costfunctional}
f(\nu) = \sum_{i=1}^L W_2^2(\nu,\nu_i)\,.
\end{equation}
The desired permutation matrices are now given by
\[
P^i = M T^i, \qquad i\in \{1,\ldots,L\}\,,
\]
where $T^i$ denotes the optimal coupling matrix associated to $W_2^2(\nu,\nu_i)$. 
Efficient numerical methods for solving Wasserstein barycenter problems for empirical measures 
have been discussed in \cite{CuturiDoucet2013}. These permutation matrices $P^i$ are now used to
rearrange the state samples prior to the application of (\ref{transform_x}). 

While the transformation steps (\ref{transform_lambda})--(\ref{transform_x}) in the parameters and the state samples is relatively complex, we emphasize that it only needs to be conducted whenever the effective sample size (\ref{effective_SS}) drops below a certain threshold such as 
$L^\ast = 3L/4$, for example.

%In fact it is possible to determine a linear transformation that maximizes the correlation via solving an optimal transport problem
%\begin{equation}\label{Mongeproblem}
%T^*=\underset{T\in\mathcal{T}(\Lambda,\tilde \Lambda)}{\arg\inf}\mathbb{E}[||\Lambda-\tilde \Lambda||^2]
%\end{equation}
%with $T\in\mathcal{T}(\Lambda,\tilde \Lambda)$ being the set of maps $T:\Lambda\rightarrow \tilde \Lambda$  such that $\tilde \pi$ is the push forward with respect to T of $\pi$ ,i.e., $T_{\#}\pi=\tilde \pi$.
%Note that the solution of the optimization problem in (\ref{Mongeproblem}) is a deterministic coupling that also induces a coupling that minimizes the Monge Kantorovich problem. In order to numerical approach (\ref{Mongeproblem}) for the filtering problem at hand the discretized version 
%\begin{equation}\label{lineartransportproblem}
%T^*=\underset{D\in\mathbb{R}^L\times\mathbb{R}^L}{\arg\min}\sum_{i,k=1}^M t_{ik} \|{\lambda}^i - {\lambda}^k\|^2
%\end{equation}
%subject to constraints
%\begin{equation}\label{constraintslineartransport}
%\sum_{i=1}^L t_{ik}=1/L\quad, \sum_{k=1}^L t_{ik}=w^i\text{ and }t_{ik}\ge 0
%\end{equation}
%has to be solved. In the following section we will discuss the algorithmic implementation of the proposed approach.

%%%%%%%%%%%%%%%%%%%%%%%%%
%
\section{Algorithmic summary} \label{sec:AS}
%
%%%%%%%%%%%%%%%%%%%%%%%%
The details of the proposed hybrid mixture model are laid out in form of pseudocode in Algorithm (\ref{alg:Mixture}). In particular, the states $\hat x_t^{i,j}$ are 
evolved numerically via a forward Euler discretization of the EnKBF (\ref{ENKBF}) 
with step-size $\Delta t$  and the weight update formula (\ref{eq:weights}) is discretized as
\begin{equation}\label{eq:disweights}
w_{n+1}^i \propto 
w_n^i\exp^{-\frac{1}{2} (\bar h^i_{n})^{\top}R^{-1}\bar h^i_{n}\Delta t-(\bar h^i_n)^{\top}R^{-1} \Delta y_n}
\end{equation}
in order to prevent negative weights. Here subscript $n$ denotes approximations at time-level $t_n = n\Delta t$. 
 
Whenever the effective sample size (\ref{effective_SS}) drops below a threshold value $L^\ast <L$, the parameters and the state samples are updated via the extended ETPF as described in
Section \ref{Sec:ETPF} and the weights are reset to $1/L$. This step requires to solve a 
linear transport problem, {\rm i.e.},
\begin{equation}\label{lineartransportproblem}
T^*=\underset{T\in\mathbb{R}^{L \times L}}{\arg\min}\sum_{i,k=1}^L t_{ik} \|
\hat{\lambda}_n^i - \hat{\lambda}_n^k\|^2
\end{equation}
subject to the constraints
\begin{equation}\label{constraintslineartransport}
\sum_{i=1}^L t_{ik}=1/L\,,\,\,  \sum_{k=1}^L t_{ik}=w_n^i\,\,\text{ and }\,\,t_{ik}\ge 0
\end{equation}
and, if necessary, the Wasserstein barycenter problem (\ref{costfunctional}).  

\begin{algorithm}[h]
 \caption{Two step EnKF-ETPF update}
 \begin{algorithmic}[1]
     \REQUIRE $\hat x_0^{i,j}$, $\hat 
     \lambda_{0}^i$, $w^i_0:=1/L$, $i\in\{1,\dots,L\}$ $j\in\{1,\dots,M\}$, $\Delta y_{[0,t]}$
     \ENSURE $\bar\lambda_{[0,t]}$, $\bar x^i_{[0,t]}$, $i\in\{1,\ldots ,L\}$ \medskip

     \medskip
     \FOR {$n=1$ to  $t$}
         \medskip
     \FOR {$j=1$ to  $M$}
       \STATE Calculate $\hat x_t^{i,j}$ according to eq. (\ref{ENKBF})
     \ENDFOR
       \medskip

     \FOR {$i=1$ to  $L$}
           \STATE Compute weights weights $w^i_n$ via (\ref{eq:disweights})\medskip
 \IF{$L_t^{\rm eff}\le L^*$}
     \STATE Solve minimization linear transport problem (\ref{lineartransportproblem}) to find $T^*$ \medskip
     \STATE If required, solve the Wasserstein Barycenter problem (\ref{costfunctional}) and rearrange
     states $\hat x_t^{i,j}$
     \medskip
       \STATE Update samples $\hat x_t^{i,j}$ and $\hat \lambda^i_n$ via (\ref{transform_x}) and (\ref{transform_lambda}) \medskip
        \STATE Set weights $w_n^i:=\frac{1}{L}$\medskip
  \ENDIF 
   \ENDFOR
     \STATE Determine $\bar \lambda_n=\frac{1}{L}\sum_{i=1}^L \hat \lambda^i_n$ $\forall~n$
     \STATE Compute $\bar x_n^i= \frac{1}{M}\sum_{j=1}^M \hat x^{i,j}_n$ $\forall~n$ $\forall~i$
        \ENDFOR

%      \[
%      \K^i = \frac{1}{2\epsilon}\sum_{j=1}^N \left[T_{ij}(\Phi_j + \epsilon(h(X^j) 
%      - \hat{h}^{(N)}))\left(X^j - \sum_{k=1}^N T_{ik}X^k\right)\right]
%      \]
 \end{algorithmic}
 \label{alg:Mixture}
\end{algorithm}

%%%%%%%%%%%%%%%%%%%%%%%%%%%%%%%%%%%%%%
%
\section{Numerical example} \label{sec:NE}

\begin{figure}
\centering
\includegraphics[width=0.8\columnwidth]{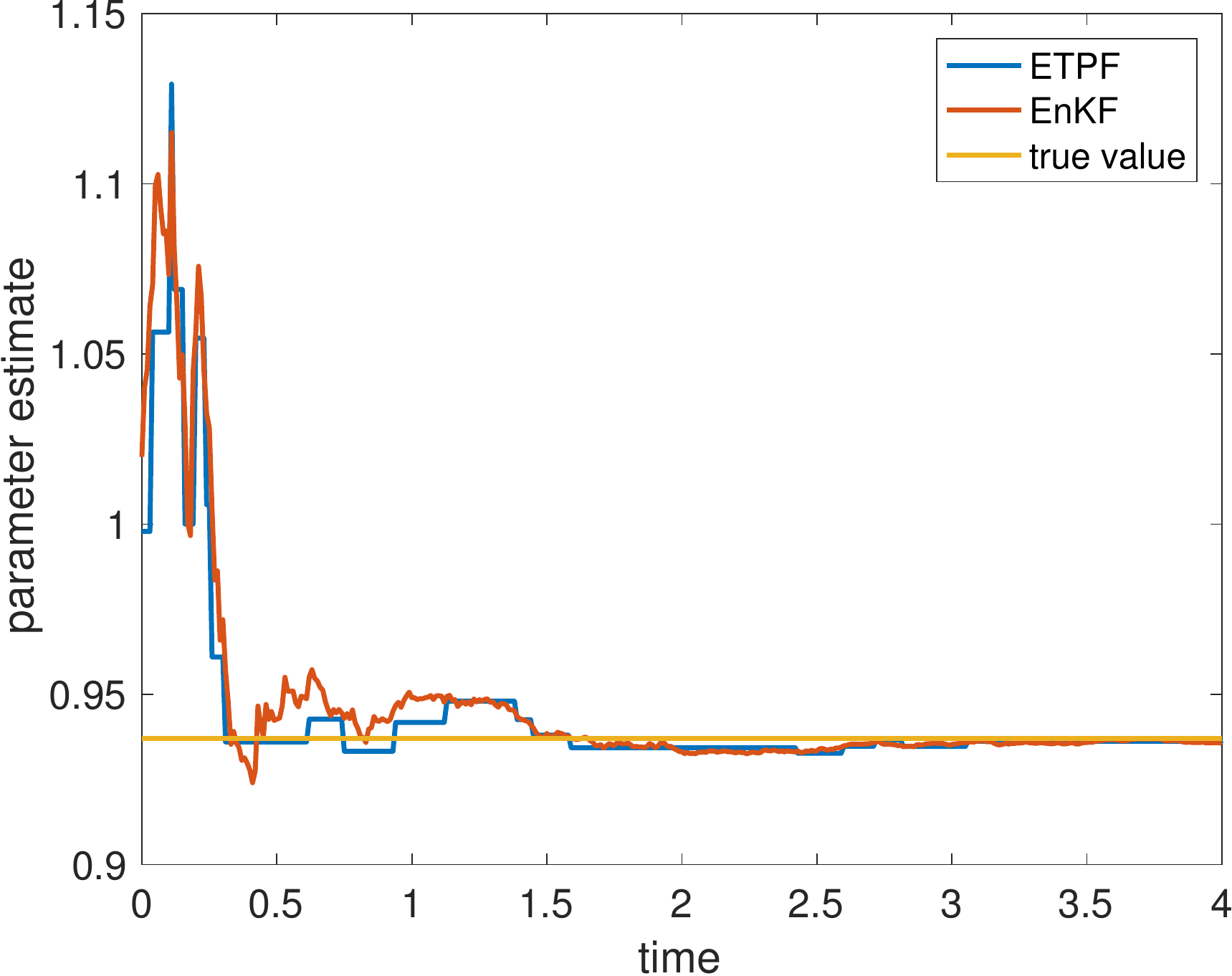}
%\vspace*{-0.6in}
\caption{Time evolution of the estimated  wave velocity parameter, $c$, is displayed for two different techniques for updating the parameter values, namely the EnKBF (red) and the ETPF (blue). 
The true parameter value  is shown in yellow. In both cases, the EnKFB is used for updating the ensembles of states.}
\vspace*{-0.1in}
\label{fig:Timeevolutionc}
\end{figure}
The proposed mixture ansatz is now numerically investigated for the stochastic wave equation 
\begin{align}\label{waveeq1}
\dot{v}_{t}&=c\bigtriangleup u_{t}+\gamma \bigtriangleup v_{t}+\delta \dot{W}_t(x)\,, \\
\dot{u}_t&=v\,, \label{waveeq2}
\end{align}
with unknown wave velocity parameter $c = e^\lambda$, $\delta=0.02$, $\gamma = 0.001$ and 
space-time white noise $\dot{W}_t(x)$. Note that a more general case than discussed earlier is considered here, {\rm i.e.},  
$u_t(x)$ and $v_t(x)$ are functions with respect to the spatial domain $x \in [0,2\pi]$ at any given time $t$. We assume periodic boundary conditions and generate initial fields from a Gaussian process prior.
 The spatial domain is discretized by means of a set of  $100$ equidistant points and the resulting
finite-dimensional SDE is integrated over the time interval $ t\in [0,4]$ with step-size $\Delta t=0.01$ ({\rm i.e.}, 400 time steps).
We observe the velocity field $v_t(x)$ continuously in time
\begin{equation}
\dot{y}_t=v_t +R^{1/2} \dot{V}_t\,,
\end{equation}
where $\dot{V}_t$ denotes space-time white noise and $R=0.0001$. The threshold for the effective sample size is set to $L^\ast = 3L/4$. 
The ensemble size  of states is $M=100$ for each parameter $c_i = e^{\hat \lambda_0^i}$ 
and we use $L=20$ different realizations of the unknown parameter $\lambda$, which is assumed to
Gaussian distributed at initial time.  
As shown in Figure \ref{fig:Timeevolutionc}, the proposed methodology 
yields a correct estimate for the true parameter value (displayed in yellow) after a relatively short assimilation window. 
We also display an estimate using an EnKFB for both the state and parameter updates for comparison. Since all $L=20$ different models where initialized with the same ensemble of states and
the data assimilation window is relatively short, the solution of a Wasserstein barycenter problem was not necessary in order to keep the different ensembles of states sufficiently correlated. 
It is not surprising that the EnKBF performs well for this state-parameter estimation problem since the conditional density $\tilde \pi_t$ is close to Gaussian. The ETPF parameter update should provide
more appropriate for non-Gaussian distributions in the parameter values.

%%%%%%%%%%%%%%%%%%%%%%%%%%%%%%%%%%%%%%

\section{Conclusions}\label{Sec:Conclusions} 

We have presented a sequential state-parameter estimation algorithm suitable for high-dimensional state space models which depend on a relatively small number of 
parameters. In comparison to a direct approach based on the FPF formulation (\ref{FPF}) for an extended state space model, the proposed methodology can be implemented
as a parallel update of $L$ standard state estimation problems for given parameters using either 
the FPF, the EnKBF or other ensemble transform particle filter. 
If necessary, techniques such localization and ensemble inflation
\cite{sr:evensen,Reich2015} can also be used. The parameters, on the other hand, are adjusted using an extended version of the ETPF once the effective sample size of the parameters 
drops below a certain threshold value. This part of the algorithm is computationally more demanding than a standard resampling approach. However, it allows again for an application 
of localization to the estimation of spatially dependent parameters. A practical exploration of such an extension of the presented algorithm will be explored for stochastic wave equation 
(\ref{waveeq1})--(\ref{waveeq2})
and spatially dependent wave velocities $c(x)$ in future work. Furthermore, the computational complexity of the parameter update step can be reduced by using alternative implementations 
of the ETPF such as provided by the Sinkhorn algorithm for the underlying linear transport problem \cite{jdw:DeWiljesAcevedoReich2017}.
%%%%%%%%%%%%%%%%%%%%%%%%%%%%%%%%%%%%%

\section*{Acknowledgments} This research has been funded by the Deutsche Forschungsgemeinschaft (DFG) through grant CRC 1294 \textit{Data Assimilation}, Project (A02) ``Long-time stability and accuracy of ensemble transform filter algorithms'' and Project (A06) ``Approximative Bayesian inference and model selection for stochastic differential equations".

%%%%%%%%%%%%%%%%%%%%%%%%%%%%%%%%%%%%%%%%%%%%%%%%%%%%%%%%%%%%%%%%%%%%

\bibliographystyle{siam}
\bibliography{Ref}

\end{document}